\documentclass[10pt,leqno]{amsart}
\usepackage{graphicx}
\usepackage{indentfirst,csquotes}

\usepackage{amssymb,amsthm,amsmath}
\usepackage{xcolor,paralist,hyperref,titlesec,fancyhdr,etoolbox}

\titleformat{\section}[display]{\normalfont\huge\bfseries\centering}{\centering\chaptertitlename\thechapter}{10pt}{\Large}
\titlespacing*{\section}{0pt}{0ex}{0ex}

\hypersetup{ colorlinks=true, linkcolor=black, filecolor=black, urlcolor=black }

\usepackage{lipsum}

\begin{document}
\title{Inverse problem of determining the right-hand side of a one-dimensional fractional diffusion  equation with variable coefficients} 
\author[D.K.Durdiev]{D.K.Durdiev}
\date{\today}
\address{$ ^1$Bukhara Branch of Romanovskii Institute of Mathematics,
Uzbekistan Academy of Sciences, Bukhara, Uzbekistan,\\
$ ^2$Bukhara State University, Bukhara, 705018 Uzbekistan}
\email{d.durdiev@mathinst.uz, durdiev65@mail.ru}
\maketitle

\let\thefootnote\relax

\begin{abstract}
In this paper, we study the inverse problem of finding a time-dependent multiplier of the right-hand side of a time-fractional  one-dimensional diffusion equation with variables coefficients in the case where the usual Cauchy, homogeneous Dirichlet boundary, and an integral overdetermination conditions are given. The overdetermination condition has the form of an integral with a weight over a spatial segment from the solution of the direct problem, in which the weight function is a spatially dependent known factor of the right-hand side of the equation. This made it possible to construct a solution to the inverse problem in explicit form and prove its correctness in the class of regular solutions
\end{abstract} 

\bigskip


\begin{center}
    {\bf Formulation of the problem}
\end{center}

In the domain $\Omega_T:=\left\{(t, x): 0< t \leq T, \, 0<x < l\right\},$ $T>0, \, l>0$ are some fixed numbers, we consider the initial-boundary value problem of determining the function $u(t, x),$ satisfying the equation
$$
  \partial_t^{\alpha}u-\frac{\partial }{\partial x}\left(p(x)\frac{\partial u }{\partial x}\right)+q(x)u=F(t, x) \eqno(1)
$$
with initial
$$
u(0, x)=\varphi(x), \, \, \, 0\leq x \leq l  \eqno(2)
$$
and boundary
$$
  u(t, 0)=u(t, l)=0, \, \, \,  0\leq t \leq T \eqno(3)
$$
conditions, where the Caputo fractional differential operator
$\partial_t^{\alpha}$ of the order $0<\alpha<1$  is defined by  [[1], pp. 90-99]:
$$
\partial_t^{\alpha}u(x, t):=I_{0+}^{1-\alpha}u_t(x, t)=\frac{1}{\Gamma(1-\alpha)}\int_0^t\frac{u_{\tau}(x, \tau)}{(t-\tau)^{\alpha}}d\tau,
$$
$$
I_{0+}^{\alpha}u(x, t):=\frac{1}{\Gamma(\alpha)}\int_0^t\frac{u(x, \tau)}{(t-\tau)^{1-\alpha}}d\tau,
$$
$I_{0+}^{\alpha}u(x, t)$  is the Riemann--Liouville fractional
integral of the  function $u(x, t)$ with respect to $t$ and $f(x,t),$ $\varphi(x)$ are given smooth functions.

Assume that the coefficients of equation (1) satisfy   the following conditions:

$A1:$ $0<p(x)\in C^1[0, l],$ \, \, \, $0\leq q(x)\in C[0, l].$

Let $C(D),$ $C\left(\overline{D}\right)$ be spaces of continuous functions on $D,$ $\overline{D},$ respectively.
Let us introduce into consideration the class of functions $C^{\alpha,\, 2}_{t, x}(D):=\Big\{v(t, x)\in C(D): \, \partial_t^{\alpha}u, \, \, u_{xx}\in  C(D) \Big\}.$

We  call (1)-(3) as a {\bf direct problem}  and we will look for its solution in the class $C^{\alpha,\, 2}_{t, x}(D)\cap C\left(\overline{D}\right).$

Assuming $F(t, x)=f(t)h(x)$, we study the following inverse problem:

\textbf{Inverse problem.} Find a  function $f(t)\in C[0, T]$ by the given overdetermination  condition for solution of the direct problem
 (1)-(3)
$$
\int\limits_0^l h(x)u(t, x)dx=g(t), \, \, \, \, 0\leq t \leq T, \eqno(4)
$$
where $ h(x)$ and $ g(t)$ are given  functions.

The last two decades, the time-fractional diffusion equations with variables coefficients have received
much attention to the most of the research communities due to their applications
in various fields of science and engineering, such as physics, chemistry, biology,
control theory, viscoelasticity and so on (see [2]-[4] and references therein).
In the case $\alpha = 1,$  problems of determining the right-hand side factor
depending on the time variable in the parabolic equations with conditions
(2), (3) and by integral overdetermination conditions of the form (4)
were considered in  [5]. Sufficient conditions for the data of the problem
were found, under which the existence and uniqueness theorems hold. A
numerical solution algorithm was also constructed and implemented. In [6], [7],
for $\alpha = 1$ inverse problem (1)-(4) was studied in a more general formulation
with the aim of proving the existence and uniqueness of a generalized solution.
Many articles were addressed for inverse problems of determination of a
time-dependent right hand side in time-fractional diffusion equation (1)
with $p(x) = 1, \, q(x) = 0$ and various type of time fractional derivative. In
this paper [8], an inverse problem of determining a time-dependent source
term in a one-dimensional time-fractional diffusion equation from the energy
measurement is studied. This paper [9] is devoted to identifying a time-
dependent source term in a multi-dimensional time-fractional diffusion equation
from boundary Cauchy data. In [10], the inverse problem of finding the time-
dependent coefficient of a generalized time fractional diffusion equation, in
the case of non-local boundary and integral overdetermination conditions
is studied. The existence and uniqueness of the solution of the considered
inverse problem are obtained by a method based on the expansion of the
solution by using a bi-orthogonal system of functions and the fractional
calculus. The paper [11] is devoted to forward and inverse source problems for a 2D in space variables time-fractional diffusion equation.
For the inverse problem for determining a time-dependence function, Abel's integral
equation of the first kind is obtained, which is further reduced to an integral equation of the second kind with the application of fractional
differentiation. Solvability results are obtained.

In this work, in contrast to the above-mentioned works, consideration of
the overdetermination integral condition (4)  $h(t)$  is the spatial component
of the acting load $F (x, t)$) made it possible to construct solutions to inverse
problem (1)-(4)  in explicit form and prove their correctness in the class
of regular solutions.
The rest of the manuscript is organized as follows: In Section 2, we investigate direct problem  based on the fractional calculus and Fourier analysis method.  Section 3 studies inverse problem and  the
existence and uniqueness of the solution to this  problem  are
obtained by reducing the problem to an integral equation of the second kind of the Abelian  type.
Section 4  contains some conclusions.

\begin{center}
    {\bf Study of the direct problem}
\end{center}

We will seek a solution to the problem (1)-(3) in the form of a Fourier series
$$
u(x, t)=\sum_{n=1}^{\infty} u_n(t) X_n(x) \eqno(5)
$$
by the eigenfunctions of the equation
$$
L(X)\equiv\frac{d }{d x}\left(p(x)\frac{d X }{d x}\right)-q(x)X=-\lambda X \eqno(6)
$$
with the boundary conditions
$$
 X(0)=X(l)=0.  \eqno(7)
$$

If conditions of $A1$ are satisfied, it is shown in books [[12], pp. 168--215],  [[13], pp. 225--238] that the spectral problem (6), (7) in $L^2(0, l)$ has a complete system of orthonormal eigenfunctions $X_n(x),$ and the positivity of the corresponding eigenvalues $\lambda_n$ can be proved by the following standard technique. For $X=X_n(x)$ and $\lambda=\lambda_n$, we scalar multiply equation (6) by $X_n(x)$ and, taking into account the normalization of the eigenfunctions, we have
$$
\lambda_n=-\int\limits_0^l\left[\frac{d }{d x}\left(p(x)\frac{d X_n }{d x}\right)-q(x)X_n\right]X_n(x)dx.
$$
Integrating the first term on the right-hand side of this equality and taking into account (7), we obtain the equality
$$
\lambda_n=\int\limits_0^l\left[p(x)\left(\frac{d X_n }{d x}\right)^2+q(x)X_n^2\right]dx, \eqno(8)
$$
which means $\lambda_n>0, \, \, n\in \mathbb{N}. $
It is clear that $\lambda_n$ form a countable set and with a suitable numbering, we have $0<\lambda_1<\lambda_2<\cdots$, $\lim\limits_{n\rightarrow\infty}\lambda_n=\infty$. In  [[12], p. 182], in particular, the estimate was obtained
$$
c_1n^2+c_2\leq \lambda_n \leq C_1n^2+C_2, \eqno(9)
$$
where pairs of constants ${c_1, c_2}$ and ${C_1, C_2}$ depend, respectively, on the minimum and maximum values of the functions $p(x), \, q(x)$ on the interval $[0, l].$

Expanding the right-hand side of (1) in terms of these functions and equating the expressions for the same eigenfunctions of the equation (1), we find ordinary fractional differential equations for determining the Fourier coefficients $u_n(t):$
$$
\partial_t^{\alpha}u_n(t)+\lambda_n u_n(t)=F_n(t), \eqno(10)
$$
where $F_n(t)= \int_0^l F(x, t) X_n(x) d x $ and $L(X_n)=-\lambda_n X_n.$ It is easy to verify that the solution of equation (10) satisfying the condition
$$
u_n(0)=\varphi_n=\int\limits_0^l \varphi(x) X_n(x) d x,
$$
for each $n\in \mathbb{N},$ is the function [[1], pp. 140, 141]
$$
u_n(t)=\varphi_n E_{\alpha, \, 1}\left(-\lambda_n t^{\alpha}\right)+\int\limits_0^t F_n(s)(t-s)^{\alpha-1} E_{\alpha, \, \alpha}\left(-\lambda_n (t-s)^{\alpha}\right) d s, \eqno(11)
$$
where   $E_{\alpha,
\beta}(z)$ is the Mittag-Leffler function defined by the following series  [[1], pp. 40-45]:
$$
E_{\alpha,
\beta}(z):=\sum_{n=0}^{\infty}\frac{z^n}{\Gamma(\alpha{n}+\beta)},
$$
 where $\Gamma(\cdot)$ is the  Euler's gamma function $\alpha,z, \beta\in\mathbb{C}$, $
\mathfrak{R}(\alpha)>0.$
This function has the following properties [[1], pp. 40-45]:

{\bf Proposition 1.}   For $0<\alpha<1$, $t>0$, we have $0<E_{\alpha,1}(-t)<1$. Moreover, $E_{\alpha,1}(-t)$ is completely monotonic, that is
$$
(-1)^n\frac{d^n}{dt^n}E_{\alpha,1}(-t)\geq0,\quad \mbox{for} \ n\in\mathbb{N}.
$$

{\bf Proposition 2.}  For $0<\alpha<1$, $\eta>0$, we have $0\leq E_{\alpha,\alpha}(-\eta)\leq\frac{1}{\Gamma(\alpha)}$. Moreover, $E_{\alpha,\alpha}(-\eta)$ is a monotonic decreasing function with $\eta>0$.

{\bf Proposition 3.} Let $0<\alpha<2$ and $\beta\in\mathbb{R}$ be arbitrary. We suppose that $\kappa$ is such that $\pi\alpha/2<\kappa<\min\{\pi,\pi\alpha\}$. Then there exists a constant $C=C(\alpha,\beta,\kappa)>0$ such that
$$
\left|E_{\alpha,\beta}(z)\right|\leq\frac{C}{1+|z|},\quad \kappa\leq|\mbox{arg}(z)|\leq\pi.
$$

Thus, the solution of the problem (1)-(3) must be expressed as a formal series
$$
u(x, t)=\sum_{n=1}^{\infty} \Big[\varphi_n E_{\alpha, \, 1}\left(-\lambda_n t^{\alpha}\right)\quad\quad\quad\quad\quad\quad\quad\quad\quad\quad\quad\quad\quad\quad\quad$$
$$
\quad\quad\quad\quad\quad\quad+\int\limits_0^t F_n(s) (t-s)^{\alpha-1} E_{\alpha, \, \alpha}\left(-\lambda_n (t-s)^{\alpha}\right) d s\Big] X_n(x). \eqno(12)
$$
If the series (12) converges uniformly in $\overline{D}$ and the series obtained from (12) by termwise applying  the fractional operator $\partial_t^{\alpha}$ and twice differentiation  with respect to $x$ converge uniformly in $D$, then the sum of this series is from the class $ C\left(\overline{D}\right)\cap C^{\alpha, 2}_{t, x}\left(D\right),$ satisfying equation (1) and conditions (2), (3).

The following statement is true.

\textbf{Theorem 1.} Let conditions $A1,$ $F(t, x)=f(t)h(x)$ and $\varphi(x) \in C^1[0, l],$ $ \, \varphi(0)=\varphi(l)=0,$ $ \, h(x) \in C^1[0, l], \, h(0)=h(l)=0,$ $ \, f(t) \in C[0, T]$ be satisfied, then there exists a unique solution to problem (1)-(3), defined by formula (12), where $u_n(t)$ is found by formula (11).

\textbf{Proof.} Each term of the series (12) satisfies the equation (1) by its construction, and it suffices for us to show the uniform and absolute convergence in a closed domain  $D_{t_0}:=\big\{(t, x): $ $ t_0\leq t \leq T, \, 0\leq x \leq l\big\},$ $t_0$ is a sufficiently small positive number, of the series
$$
u(t, x)=\sum\limits_{n=1}^{\infty} u_n(t) X_n(x), \, \, \partial_t^{\alpha}(x, t)=\sum\limits_{n=1}^{\infty} \partial_t^{\alpha}u_n(t) X_n(x), \eqno(13)
$$
$$
u_x(t, x)=\sum\limits_{n=1}^{\infty} u_n(t) X'_n(x), \, \, u_{xx}(x, t)=\sum\limits_{n=1}^{\infty} u_n(t) X''_n(x),  \eqno(14)
$$
where $\partial_t^{\alpha}u_n(t)$ as follows from (10)
$$
\partial_t^{\alpha}u_n(t)=-\lambda_n u_n(t)+h_nf(t), \, \, h_n= \int_0^l h(x) X_n(x) d x.   \eqno(15)
$$

  Since $F_n(t)=h_n\,f(t)$ it follows from (12) based on  Propositions 1 and 2  that the first series of  (13) in $D_{t_0}$ are less in absolute value than the expression
  $$
  C_0\left(\sum_{n=1}^{\infty}\left|\varphi_n X_n(x)\right|+\sum_{n=1}^{\infty}\left|h_n X_n(x)\right|\right), \eqno(16)
  $$
  where the constant $C_0$ depends on $T,$ $\alpha$ and $f(t).$

To prove the convergence, for example, of the first series in (16), first we carry out the following calculations.

Let us consider the functional, defining it by the  right-hand side (8):
$$
J(Y)=\int\limits_0^l\left[p(x) \left(Y'(x)\right)^2+q(x)Y^2(x)\right]dx.
$$
Substituting into it the expression $Y(x)=\varphi(x)-\sum\limits_{j=1}^{n-1}\varphi_j X_j(x) $ and taking into account (8), as well as the easily verifiable identity
$$
    \int\limits_0^l\left[p( x) X'_k(x) X'_j(x)+q(x)X_k(x)X_j(x)\right]dx=0, \, \, \, k\neq j,
$$
we get
$$
J\left[\varphi(x)-\sum\limits_{j=1}^{n-1}\varphi_jX_j(x)\right]=\int\limits_0^l\left[p(x)\left(\varphi'(x)\right)^2+q(x)\varphi^2(x)\right]dx
$$
$$
+\sum\limits_{j=1}^{n-1}\lambda_j\varphi_j^2-2\sum\limits_{j=1}^{n-1}\varphi_j\int\limits_0^l\left[p(x)\varphi'(x) X'_j(x)+q(x)\varphi(x)X_j(x)\right]dx.
$$
Given the conditions $\varphi(0)=\varphi(l)=0,$ we integrate by parts in the last integral. Taking into account the equality (6) for $X=X_n(x),$  $\lambda=\lambda_n,$   we find
$$
J\left[\varphi(x)-\sum\limits_{j=1}^{n-1}\varphi_jX_j(x)\right]=\int\limits_0^l\left[p(x)\left(\varphi'(x)\right)^2+q(x)\varphi^2(x)\right]dx
-\sum\limits_{j=1}^{n-1}\lambda_j\varphi_j^2. \eqno(17)
$$
Since $J(Y)\geq 0,$ then from formula (17) the inequality
$$
\sum\limits_{j=1}^{\infty}\lambda_j\varphi_j^2\leq\int\limits_0^l\left[p(x)\left(\varphi'(x)\right)^2+q(x)\varphi^2(x)\right]dx, \eqno(18)
$$
and the convergence of the series on the left-hand side follow directly.

Next, using the integral equation for the eigenfunctions $X_n(x)$
$$
X_n(x)=\lambda_n\int\limits_0^lG(x, \xi)X_n(\xi) d\xi,   \eqno(19)
 $$
 where $G(x, \xi)$ is the Green's function of the problem (6) and (7), we represent the series \\ $ \sum_{n=1}^{\infty}\left|\varphi_n X_n(x)\right|$ in the form
$$
\sum_{n=1}^{\infty}\lambda_n\left|\varphi_n \gamma_n(x)\right|, \,\,\,\,\,  \gamma_n(x)=\int\limits_0^lG(x, \xi)X_n(\xi)\,d\xi,   \eqno(20)
$$
$\gamma_n(x)$  are the  Fourier coefficients of the function $G(x, \xi)$ with respect to the argument $\xi.$ For them, according to inequality (18), the formula
$$
\sum\limits_{j=1}^{\infty}\lambda_j\gamma_j^2\leq\int\limits_0^l\left[p(x)G^2_{\xi}(x, \xi)+q(x)G^2(x, \xi)\right]dx.
$$
Hence, due to the boundedness of the integrands, it follows the uniform boundedness of the series
$$
\sum\limits_{j=1}^{\infty}\lambda_j\gamma_j^2\leq M=const. \eqno(21)
$$
Then, applying the Cauchy-Bunyakovsky inequality to the segment of the series (20), we have
$$
\sum_{n=m}^{m+k}\lambda_n\left|\varphi_n \gamma_n(x)\right|=\sum_{n=m}^{m+k}\sqrt{\lambda_n}\left|\varphi_n\right|\sqrt{\lambda_n}\left| \gamma_n(x)\right| \leq \sqrt{\sum_{n=m}^{m+k}\lambda_n\varphi^2_n}\,\sqrt{\sum_{n=m}^{m+k}\lambda_n \gamma^2_n(x)}
$$ or
$$
\sum_{n=m }^{m+k}\lambda_n\left|\varphi_n \gamma_n(x)\right| \leq \sqrt{\sum_{n=m}^{m+k}\lambda_n\varphi^2_n}\,\sqrt{M}.
$$
From this inequality and the convergence of the series with terms $\lambda_n\varphi^2_n$ it follows immediately that the series (20) converges uniformly for $x\in[0, l],$ i.e. the first series in (16) converges uniformly. If the conditions of Theorem 1 are satisfied with respect to the function $h(x),$ the proof of the convergence of the second series in (16) literally repeats the calculations performed.

Further it follows from (12) based on  Propositions 3  that the second  series of  (13) in $D_{t_0}$ are less in absolute value than the following  expression
$$
C_1\left(\sum_{n=1}^{\infty}\frac{\lambda_n}{1+\lambda_n t_0^{\alpha}}\left|\varphi_n X_n(x)\right|+\sum_{n=1}^{\infty}\frac{\lambda_n}{1+\lambda_n t_0^{\alpha}}\left|h_n X_n(x)\right|+\sum_{n=1}^{\infty}\left|h_n X_n(x)\right|\right), \eqno(22)
$$
where the constant $C_1$ depends on $T,$ $\alpha$ and $f(t).$ Since $\frac{\lambda_n}{1+\lambda_n t^{\alpha}}<t_0^{-\alpha},$ then in this case the series (22) is estimated by the series  (16) (with another constant instead of $C_0$), for which the convergence were shown.

First, note that the Green's function $G(x, \xi) $ has the form
$$
G(x, \xi)=y_1(x)\,y_2(\xi) \, \, \, \text{at} \, \, \, 0 \leq \xi \leq x; \, \, \, \, G(x, \xi)=y_2(x)\,y_1(\xi) \, \, \, \text{at} \, \, \, x\leq \xi \ leq l,
$$
where $y_1(x)$ and $y_2(x)$ are linearly independent partial solutions of the problem (6), (7).
To study the series (14), using this, we write the expression (19) for $X_n(x)$ in the form
$$
X_n(x)=\lambda_n y_1(x)\int\limits_0^xy_2(\xi)X_n(\xi)d\xi+\lambda_n y_2(x)\int\limits_0^ly_1(\xi)X_n(\xi)d\xi.
$$
From here, by direct differentiation and then multiplication by $\varphi_n$, due to the equality $\left(p(x)X_n'(x)\right)'-q(x)X_n(x)=-\lambda_n X_n(x)$, we obtain
$$
\varphi_nX'_n(x)
= y'_1(x)\int\limits_0^xy_2(\xi)\lambda_n \varphi_n X_n(\xi)d\xi+ y'_2(x)\int\limits_0^ly_1(\xi)\lambda_n \varphi_nX_n(\xi)d\xi. \eqno(23)
$$

$$
\varphi_nX''_n(x)=\frac{1}{p(x)} \Big[-p'(x)X'_n(x)\varphi_n + q(x)\varphi_nX_n(x)-\lambda_n\varphi_nX_n(x)\Big]. \eqno(24)
$$
It is clear that the series (14) in $D_{t_0}$ in absolute value do not exceed the following expressions, respectively:
$$
C_2\left(\sum_{n=1}^{\infty}\left|\varphi_n X'_n(x)\right|+\sum_{n=1}^{\infty}\left|h_n X'_n(x)\right|\right), \eqno(25)
$$
$$
C_3\left(\sum_{n=1}^{\infty}\left|\varphi_n X''_n(x)\right|+\sum_{n=1}^{\infty}\left|h_n X''_n(x)\right|\right), \eqno(26)
$$
where constants $C_i,$ $i=2,3,$ depend on $T,$ $\alpha$ and $f(t).$

Taking into account uniform convergence of the second series in (13), and also substituting (23) into the first series in (14) and representation (23) with $\varphi_n$ replaced by $h_n$ into the second series in (14), we verify that the series in (25) converge uniformly. From equality (24)
and convergence of the series in (25) follows uniform convergence in $D_{t_0}$ of the series in (26). It remains to note that the series (12) converges absolutely and uniformly for $t\rightarrow 0$. Its sum is a continuous function for $x\in[0, l],$ i.e.
$$
\lim\limits_{t\rightarrow 0}u(t, x)=u(0, x)=\sum\limits_{n=1}^{\infty} \varphi_n X_n(x)=\varphi(x).
$$
Theorem 1 is proved.

\begin{center}
    { \bf Investigation of the inverse problem}
\end{center}

Let $AC^n[0, T]$ be  the class of continuously differentiable  on $[0, T]$ up to the order of $n-1$ functions $f(x)$ and $f(x)^{(n-1)}\in AC[0, T],$ where $n=1,2, ...$ and $AC[0, T]$ is the class of absolutely continuous   on $[0, T]$ functions.

The following statement is true.

{\bf Theorem 2.} Let the functions $\varphi(x)$ and $h(x)$ satisfy the conditions of Theorem 1, $h(x) \not \equiv$ 0 on $[0, l], \, g(t) \in AC^2[0, T]$ and the conditions
$$
\int_0^l \varphi(x) h(x) d x=g(0)=\sum_{n=1}^{\infty}h_n \varphi_n
$$
are valid. Then there exists a unique solution to inverse problem such that $f(t)\in AC[0, T].$

{\bf Proof.}
The solution of the direct problem (1)-(3) satisfies the additional condition (4). Then, we get
$$
\int\limits_0^l h(x)u(t, x)dx=\int\limits_0^l h(x) \sum\limits_{n=1}^{\infty} u_n(t) X_n(x) d x=\sum_{n=1}^{\infty}h_n u_n(t), \eqno(27)
$$
where $h_n=\int_0^l h(x) X_n(x) d x.$ For $F(x, t)=f(t) h(x)$ function (11) take the form
$$
 u_n(t)=\varphi_n E_{\alpha, \, 1}\left(-\lambda_n t^{\alpha}\right)+h_n f_n(t), \eqno(28)
$$
where $ f_n(t)=\int\limits_0^t f(s) (t-s)^{\alpha-1} E_{\alpha, \, \alpha}\left(-\lambda_n (t-s)^{\alpha}\right) d s.$

In (27) substituting (28) and using (4) after simple transformations we obtain the Abelian type   integral equation of the first kind with a difference kernel with respect to the unknown function $f(t)$
$$
\int\limits_0^t f(s) K(t-s) d s=G(t), \quad 0 \leqslant t \leqslant T, \eqno(29)
$$
where
$$
K(t)=t^{\alpha-1}\sum\limits_{n=1}^{\infty} h_n^2  E_{\alpha, \, \alpha}\left(-\lambda_n t^{\alpha}\right), \quad (t, s)\in \triangleleft(T)=\{0 \leqslant s \leqslant t \leqslant T\}, \eqno(30)
$$
$$
G(t)=g(t)-\sum_{n=1}^{\infty}h_n \varphi_n E_{\alpha, \, 1}\left(-\lambda_n t^{\alpha}\right), \quad 0 \leqslant t \leqslant T. \eqno(31)
$$

Due to the conditions imposed on the functions $\varphi(x)$ and $h(x)$ in Theorem 1,  the series in (30) and (31) converge uniformly and admit term-by-term differentiation with respect to $t$.

Further we use the two assertion (see, respectively, [14] and  [[15], Lemma 2.1]):

{\bf Lemma 1.}
For $0<\alpha<1$ and $\lambda>0,$ if $\eta(t)\in AC[0, T],$ then, we have
$$
\partial_t^{\alpha}\int\limits_0^t\eta(s)(t-s)^{\alpha-1}E_{\alpha, \alpha}\left(-\lambda (t-s)^{\alpha-1}\right)ds\quad\quad\quad\quad\quad\quad\quad\quad\quad\quad\quad\quad
$$$$
\quad\quad\quad\quad\quad=\eta(t)-\lambda\int\limits_0^t\eta(s) (t-s)^{\alpha-1}E_{\alpha, \alpha}\left(-\lambda (t-s)^{\alpha-1}\right)ds, \quad 0<t\leq T.
$$
In particular, if $\lambda=0,$ we have
$$ \partial_t^{\alpha}\int\limits_0^t\eta(s) (t-s)^{\alpha-1}=\Gamma(\alpha)\eta(t), \quad 0<t\leq T.$$

{\bf Lemma 2.}
For $0<\alpha<1$ , if $f(t)\in AC[a, b],$ then,
$$f_{1-\alpha}(x):=\frac{1}{\Gamma(1-\alpha)}\int\limits_a^x\frac{f(t)d t}{(t-\alpha)^{\alpha}}\in AC[a,b],$$
and
$$f_{1-\alpha}(x)=\frac{1}{\Gamma(2-\alpha)}\left[f(a)(x-a)^{1-\alpha}+\int\limits_a^xf'(t)(t-\alpha)^{1-\alpha}d t\right].$$

Applying the operator $\partial_t^{\alpha}$ to both parts of (29) and in view of  Lemma 1, we get
$$
f(t)\sum\limits_{n=1}^{\infty} h_n^2 +\int\limits_0^t \frac{K_0(t-s)}{(t-s)^{1-\alpha}}f(s)  d s=G(t), \quad 0 \leqslant t \leqslant T, \eqno(32)
$$
where
$$
K_0(t)=-\sum\limits_{n=1}^{\infty}\lambda_n h_n^2  E_{\alpha, \, \alpha}\left(-\lambda_n t^{\alpha}\right), \quad (t, s)\in \triangleleft(T)=\{0 \leqslant s \leqslant t \leqslant T\}, \eqno(33)
$$
$$
G(t)=\partial_t^{\alpha}g(t)+\sum_{n=1}^{\infty}\lambda_nh_n \varphi_n E_{\alpha, \, 1}\left(-\lambda_n t^{\alpha}\right), \quad 0 \leqslant t \leqslant T. \eqno(34)
$$
In the last formula, the equality $\partial_t^{\alpha}E_{\alpha, \, 1}\left(-\lambda_n t^{\alpha}\right)=-\lambda_nE_{\alpha, \, 1}\left(-\lambda_n t^{\alpha}\right)$ was used.
By Theorem 2, the function $h(x) \neq 0$ on $(0, l)$, therefore, according to the Parseval-Steklov equality, we have
$$
\sum_{n=1}^{\infty} h_n^2=\int_0^l h^2(x) d x=\|h\|_{L_2(0, l)}^2>0.
$$

Let us prove that the series in (33),  (34) converge absolutely and uniformly, and thus represent continuous functions.
For the series in (33), according to Proposition 2, we have
$$
\left|-\sum\limits_{n=1}^{\infty}\lambda_n h_n^2 E_{\alpha, \, 1}\left(-\lambda_n t^{\alpha}\right)\right|\leq \sum\limits_{n=1}^{\infty}\lambda_n h_n^2, \quad (t, s)\in \triangleleft(T).
$$

To prove the convergence of the  numerical series $\sum\limits_{n=1}^{\infty}\lambda_n h_n^2$, consider the following functional, defining it by the  right-hand side of (8):
$$
J(Y)=\int\limits_0^l\left[p(x)\left(\frac{d Y }{d x}\right)^2+q(x)Y^2\right]dx.
$$
Substituting into it the expression $ Y(x)=h(x)-\sum\limits_{j=1}^{n-1}h_jX_j(x) $
and taking into account (8), as well as the easily verified  identity
$$
\int\limits_0^l\left[p(x) \left(\frac{d X_k(x) }{d x}\right)\left(\frac{d X_j(x) }{d x}\right)+q(x)X_k(x)X_j(x)\right]dx=0, \, \, \, k\neq j,
$$
we get
$$
J\left[h(x)-\sum\limits_{j=1}^{n-1}h_jX_j(x)\right]=\int\limits_0^l\left[p(x)\left(\frac{d h(x) }{d x}\right)^2+q(x)h^2(x)\right]dx
$$
$$
+\sum\limits_{j=1}^{n-1}\lambda_jh_j^2+2\sum\limits_{j=1}^{n-1}h_j\int\limits_0^l\left[p(x)\frac{d h(x) }{d x}\frac{d X_j(x) }{d x}+q(x)h(x)X_j(x)\right]dx.
$$
Taking into account the conditions $h(0)=h(l)=0,$ we integrate by parts in the last integral. Taking into account the equality (6) for $X=X_n(x), $ $\lambda=\lambda_n$ and we find
$$
J\left[h(x)-\sum\limits_{j=1}^{n-1}h_jX_j(x)\right]=\int\limits_0^l\left[p(x)\left(\frac{d h(x) }{d x}\right)^2+q(x)h^2(x)\right]dx
-\sum\limits_{j=1}^{n-1}\lambda_jh_j^2. \eqno(35)
$$
Since $J(Y)\geq 0,$ then from formula (35) it follows the inequality
$$
\sum\limits_{j=1}^{\infty}\lambda_jh_j^2\leq\int\limits_0^l\left[p(x)\left(\frac{d h(x) }{d x}\right)^2+q(x)h^2(x)\right]dx,
$$
and the convergence of the series on the left.

To prove the uniform convergence of the series in (34), it is sufficient to prove the convergence of the mojarant numerical series $\sum_{n=1}^{\infty} \lambda_n |h_n| |\varphi_n|.$ Using the Cauchy inequality, we have
$$
\sum\limits_{n=1}^{\infty}\lambda_n|h_n|  |\varphi_n|=\sum\limits_{n=1}^{\infty}\sqrt{\lambda_n}|h_n|  \sqrt{\lambda_n}|\varphi_n|\leq\sqrt{\sum\limits_{n=1}^{\infty}\lambda_nh_n^2}\, \,   \sqrt{\sum\limits_{n=1}^{\infty}\lambda_n\varphi_n^2}.
$$
The convergence of the series on the right side of this relation follows directly from the previous reasoning.

Since $g(t) \in AC^2[0, T]$, then according to Lemma 2 we have $\partial_t^{\alpha}g(t) \in AC[0, T]$. The absolute continuity  of uniformly convergent series in  (33) and   (34) is obvious. Thus, the equation (32) is an integral equation of the second kind of Abelian type with absolutely continuous input data $K_0(t)$ and $G(t)$. Such an equation has a unique solution in the class of absolutely continuous functions. The solution can be found, for example, by the method of successive approximations. Theorem 2 is proved.

\begin{center}
    {\bf Conclusion}
\end{center}

 In this work, using the Fourier method, the method of integral
inequalities and the fixed point principle, the existence and uniqueness of a solution
to the inverse problem of determining the kernel of a multidimensional third-order
integro-differential pseudo-parabolic equation by an additional condition, specified
at a fixed point concerning to the solution of the first boundary value problems
proved. Apparently, all the results of this article are correct in the case when the operator $L=:\frac{\partial }{\partial x}\left(p(x)\frac{\partial  }{\partial x}\right)+q(x)$ in
(1) is replaced by the more general operator $A$, where $A$ is a self-adjoint differential
operator, defined in the domain $\Omega\subset \mathbb{R}^n$, given by $A=\sum\limits_{i, j=1}^n\frac{\partial}{\partial x_i}\left[a_{ij}(x)\frac{\partial}{\partial x_j}\right]-c(x)$
 such
that $a_{ij}(x) = a_{ji}(x),$ \, \, $\sum\limits_{i, j=1}^na_{ij}(x)\xi_i\xi_j\geq\delta\sum\limits_{i}^n\xi_i^2,$ \,\, $\delta=const> 0.$
 In addition, it is
assumed $a_{ij}(x), \, c(x)$ satisfy some conditions of smoothness and $c(x)\geq 0$. For now, this problem is open.

\end{document}